\documentclass[12pt,a4paper]{article}
\usepackage{amsthm, amscd, amsfonts, amssymb, graphicx, color}
\usepackage[pagebackref=false,colorlinks,linkcolor=red,citecolor=green]{hyperref}
\usepackage{graphicx}
\usepackage{color}
\usepackage{graphicx}
\usepackage{graphics}
\usepackage{makeidx}
\usepackage{showidx}
\usepackage{latexsym}
\usepackage{amssymb}
\usepackage{verbatim}
\usepackage{amsmath}
\usepackage{amsthm}
\usepackage{amsfonts}
\usepackage{amssymb,amsmath}
\usepackage{latexsym,amsthm,amscd}
\usepackage[all]{xy}
\newtheorem{thm}{Theorem}
\newtheorem{prop}[thm]{Proposition}
\newtheorem{defn}{Definition}
\newcounter{alphthm}
\newtheorem{Atheorem}[alphthm]{Theorem}
\newtheorem{Aproposition}[alphthm]{Proposition}

\newtheorem{cor}{Corollary}

\newtheorem{rem}{Remark}

\newtheorem{lem}[thm]{Lemma}

\newcommand{\be}{\begin{equation}}
\newcommand{\ee}{\end{equation}}
\newcommand{\ben}{\begin{enumerate}}
\newcommand{\een}{\end{enumerate}}
\newcommand{\beq}{\begin{eqnarray}}
\newcommand{\eeq}{\end{eqnarray}}
\newcommand{\beqn}{\begin{eqnarray*}}
\newcommand{\eeqn}{\end{eqnarray*}}

\newcommand{\bpf}{\begin{proof}}
\newcommand{\epf}{\end{proof}}
\newcommand{\bl}{\begin{lem}}
\newcommand{\el}{\end{lem}}
\newcommand{\bp}{\begin{prop}}
\newcommand{\ep}{\end{prop}}
\newcommand{\bd}{\begin{defn}}
\newcommand{\ed}{\end{defn}}
\newcommand{\bt}{\begin{thm}}
\newcommand{\et}{\end{thm}}
\newcommand{\bc}{\begin{cor}}
\newcommand{\ec}{\end{cor}}
\newcommand{\br}{\begin{rem}}
\newcommand{\er}{\end{rem}}
\def\nn{\nonumber}
\newcommand\bpr{\begin{prop}}
\newcommand\epr{\end{prop}}

\begin{document}
  \title{Complete Finsler spaces of constant negative Ricci curvature}
 \author{Behroz Bidabad\thanks{The corresponding author, bidabad@aut.ac.ir}, \  Maryam Sepasi }
 \date{}
\maketitle
\begin{abstract}
Here, using the projectively invariant pseudo-distance and Schwarzian derivative, it is shown that every connected complete Finsler space of the constant negative Ricci scalar is reversible.
In particular, every complete Randers metric of constant negative Ricci (or flag) curvature  is Riemannian.
\end{abstract}
\emph{Keywords:} Schwarzian; Finsler; projective; Ricci tensor; distance; reversible; Randers.
{\emph{Math. Sub. Class.} :53C60;   58B20.\\
\section{Introduction}
The Finsler metrics of negative flag curvature are widely studied by some geometres and many significant results are obtained. Without pretending to be exhaustive, we mention here the most significant ones.

Akbar-Zadeh established the following rigidity theorem in 1988;
Let $(M,F)$ be a compact connected boundary-less Finsler space of constant negative flag curvature, then $(M,F)$ is Riemannian, see \cite{Ak}.
Later, important works have been done by imposing other assumptions.
 For instance, P. Foulon has addressed the case of the strictly negative flag curvature by imposing the following additional hypothesis. The curvature is considered to be constant along a distinguished vector field on the homogeneous bundle of tangent half lines.
   He shows that under these conditions, the Finsler space is Riemannian \cite{F1}.
   Later,  he proved that a symmetric  compact Finsler space with negative curvature would be  isometric to  a locally symmetric negatively curved  Riemann space. See \cite{F2,FF} for more details.
    Next, Z. Shen considers the case of the negative but not necessarily constant flag curvature  by imposing that the $S$-curvature is constant and  proves that the Akbar-Zadeh rigidity theorem still holds in this case, see \cite{Sh}.
    In a joint work, Mo and Shen show that a compact isotropic Finsler space with negative scalar curvature is bounded above by $``-1"$ is of the Randers type \cite{MSh}.

In the present work, we have used the well-known Schwarzian derivative on Finsler manifolds, previously studied by current authors, to obtain a rigidity theorem for a complete Finsler manifold. This is a generalization, in some senses, of a result due to Akbar-Zadeh.
\bt\label{Th;main}
 Let $(M,F)$ be a connected  complete Finsler space of constant negative Ricci scalar. Then $(M,F)$ is reversible.
\et
 The ideas presented in this article are based on certain characteristics of  projectively invariant pseudo-distance defined by the well-known Schwarzian derivative operator on a Finsler manifold and are derived from our attempts in \cite{SB2,SB3} and \cite{SB1}.

 Theorem\ref{Th;main} carries out the following appealing result on Randers spaces.
\bc\label{Cor;1}
Every complete Randers metric of constant negative Ricci curvature (in particular, of constant negative flag curvature) is Riemannian.
\ec

It's noteworthy to recall that Remark\ref{Rem;2} in the page \pageref{Rem;2} shows that this result on Randers spaces has no conflict with the  one obtained by X. Cheng and Z. Shen in  \cite{chsh}.

Recall that, here our particular interest is Einstein metrics with constant negative Ricci scalar. Clearly, all the results on Einstein metrics with constant negative Ricci scalar, hold well for Finsler metrics of constant negative flag curvature.

 We close this article by showing that a Finsler metric of negative-definite parallel Ricci tensor  has reversible geodesics, see Theorem\ref{Th;2}.
\section{Preliminaries}
Let $M$ be an $n$-dimensional   connected smooth manifold. We denote by $T_xM$ the tangent space at $x\in M$, and by $TM:={\cup}_{x\in M} T_xM$ the  bundle of tangent spaces. Each element of $TM$ has the form $(x,y)$, where $x\in M$ and $y\in T_xM$. The natural projection $\pi:TM\rightarrow M $, is given by $\pi (x,y):= x$.  The pull-back tangent bundle ${\pi}^{*} TM$ is a vector bundle over the slit tangent bundle $TM_0:=TM\backslash\{0\}$ for which the fiber ${\pi}^{*}_vTM$ at $v \in TM_0$ is just $T_xM$, where $\pi (v) = x$. \\
A (globally defined) Finsler structure on $M$ is a function $F: TM\rightarrow [0 , \infty) $ with the  properties; (I) Regularity: $F$ is $C^{\infty}$ on the entire slit tangent bundle $TM_0$; (II) Positive homogeneity: $F(x , \lambda y) = \lambda F(x , y)$ for all $\lambda > 0$; (III) Strong convexity: The $n \times n$ Hessian matrix $(g_{ij}):= ({[\frac{1}{2}F^2]}_{y^iy^j})$, is positive-definite at every point of $TM_0$.
For  any $y\in T_xM_0$, the Hessian $g_{ij}(y)$ induces an inner product $g_y$ in $T_xM$ by $g_y(u,v):=g_{ij}(y)u^iv^j$.

 Let $$G^i(x,y):=\frac{1}{4}g^{ih}(x,y)(\frac{\partial^2 F^2}{\partial y^h \partial x^j} y^j -\frac{\partial F^2}{\partial x^h}) (x,y),$$
 where, $G^i$s are positively homogeneous of degree two with respect to $y^i$ and they are the coefficients of a globally defined spray $G=y^i\frac{\partial}{\partial x^i}-2 G^i\frac{\partial}{\partial y^i}$ on $TM$.
 For a vector $y \in T_xM_0$, the \emph{ Riemann curvature}
$\textbf{R}_y : T_xM \rightarrow T_xM$ is defined by $\textbf{R}_y(u)=R^i_k u^k \frac{\partial}{\partial x^i}$, where $R^i_k(y):=2\frac{\partial G^i}{\partial x^k}-\frac{{\partial}^2G^i}{\partial y^k \partial x^j}y^j+2G^j \frac{{\partial}^2G^i}{\partial y^k \partial y^j}-\frac{\partial G^i}{\partial y^j}\frac{\partial G^j}{\partial y^k}$.\\
For a two-dimensional plane $P \subset T_pM$ and a non-zero vector $y\in T_pM$, the flag curvature $\textbf{K}(P, y )$ is defined by $\textbf{K}(P,y):=\frac{g_y(u,\textbf{R}_y(u))}{g_y(y,y)g_y(u,u)-{g_y(y,u)}^2}$, where $P=$ span$\{y,u\}$. $F$ is said to be of \textit{scalar curvature} $\textbf{K}=\lambda(y)$ if for any $y\in T_pM$, the flag curvature $\textbf{K}(P,y)=\lambda(y)$ is independent of $P$ containing $y\in T_pM$. It is equivalent to the following system in a local coordinate system $(x^i,y^i)$ on                                                                               $TM$.
\begin{equation}\label{scalar cur}
R^i_k=\lambda F^2\{{\delta}^i_k-F^{-1}F_{y^k}y^i\}.
\end{equation}
If $\lambda$ is a constant, then $F$ is said to be of \textit{constant curvature}.

\subsection{A little more details on geodesics of a Finsler metric}
 By a \textit{geodesic} of a Finsler metric $F$ we understand a regular curve $\gamma: I\rightarrow M$ such that  for any sufficiently close points $a,b\in I$, $a\leq b$ the restriction of the curve $\gamma$ to the interval $[a,b]\subseteq I$ is an extremal of the length functional
 \begin{equation}\label{lenght}
 L_F(c):= \int_a^b F(c(t),\dot{c}(t)) dt,
 \end{equation}
 in the set of all smooth curves $c:=[a,b]\rightarrow M$ with $c(a)=\gamma(a)$ and $c(b)=\gamma(b)$.\\
 Note that for any orientation-preserving  re-parametrization $\tau(t)$ (that is, $\frac{d\tau}{dt}>0$), the curve $\gamma(\tau(t))$ is also a  geodesic.\\
 A geodesic $\gamma$ of a Finsler metric $F$ involves the variational problem for the length functional given by (\ref{lenght}). Therefore, it should satisfy  the well known  Euler-Lagrange equations
 \begin{equation}\label{euler}
 \frac{d}{dt} \frac{\partial F}{\partial {\dot{x}}^i}-\frac{\partial F}{\partial x^i}=0.
 \end{equation}
By direct calculations we obtain
\beq
\frac{d^2x^i}{dt^2}&+&\frac{1}{2}g^{ih}(x,y)(\frac{\partial^2 F^2}{\partial y^h\partial x^j}\frac{dx^j}{dt}-\frac{\partial F^2}{\partial x^h} (x,y))
= \frac{1}{F}(\frac{\partial F}{\partial x^j} \frac{dx^j}{dt}+ \frac{\partial F}{\partial y^j}\frac{d^2x^j}{dt^2}) \frac{dx^i}{dt}\nn\\ &&= \frac{d(log F)}{dt}\frac{dx^i}{dt}.\nn
\eeq
Finally,   geodesics are given by the following equations
\begin{equation}\label{geo*}
\frac{d^2x^i}{dt^2}+2 G^i(x(t), \dot{x}(t))=\frac{d(log F)}{dt}\frac{dx^i}{dt}.
\end{equation}
 Along a curve  $\gamma:=x^i(t)$, the \textit{arc length  parameter $``s"$} is defined by
$s(t):=  \int_a^t F(x(r),\dot{x}(r)) dr$ or $\frac{ds}{dt}=F(x(t), \dot{x}(t))$.  We have $\frac{ds}{dt}=F(x(t), \dot{x}(t))>0$ therefore it preserves the orientation. Considering  the arc length parameter $``s"$, the equation of a geodesic becomes
\begin{equation}\label{geo}
\frac{d^2x^i}{ds^2}+2 G^i(x(s), \dot{x}(s))=0.
\end{equation}
Let us consider a curve  $\gamma$ which satisfies the equations (\ref{geo}) for some parameter $``s"$. Take  a new parameter $``t"$  determined by
$$\frac{d^2s}{dt^2}=\frac{d(log F)}{dt}\frac{ds}{dt}.$$
It is clear that $\gamma$ satisfies (\ref{geo*}) and hence is a geodesic.
Now, consider a curve $\gamma$ which satisfies the equations $$\frac{d^2x^i}{dt^2}+2 G^i(x(t), \dot{x}(t))=h(t) \frac{dx^i}{dt},$$ for some differentiable real function $h$ with respect to the parameter $``t"$. Take a new parameter $``s"$ determined by $h(t)=\frac{d(log F)}{ds}.$
 According to the latter discussion, $\gamma$ is a geodesic.
It is worth noting that along a geodesic $\gamma:=(x^i(t))$ satisfying (\ref{geo}) we have $F(x(t),x^i(t))= constant$. Furthermore,  if $``s"$ is the arc length parameter we have $F(x(t),x^i(t))= 1$.

A Finsler structure $F$ is called \emph{reversible} if the opposite tangent vectors have the same length. More precisely, on a smooth manifold $M$,  a Finsler structure $F$  is called \emph{reversible} if $F(x,y)=F(x,-y)$ for all $x \in M$ and $y \in T_xM$. A Finsler structure $F$ is called with \emph{reversible geodesics}
if for any  geodesic $\gamma:t\rightarrow \gamma(t);\ t\in(a,b)$ of $F$, the reverse curve
$\bar{\gamma}: t \rightarrow \gamma(a+b-t);\  t\in (a,b)$ is also a geodesic of $F$. Clearly, every reversible Finsler structure is with reversible geodesic but the converse is not always true. For example,  every projectively flat Randers metric is geodesically reversible \cite{Shi}.
If $F$ is  reversible, then all geodesics of $F$ are reversible.
$F$ is said to be \emph{forward (resp. backward) geodesically complete}, if every  geodesic on
an open interval $(a, b)$ can be extended to a geodesic on $(a,\infty)$ (resp. $(-\infty, b)$).
$F$ is said to be \emph{complete} if it is forward and backward complete.

\subsection{Finslerian distance function}
Let  $\gamma : [a , b] \rightarrow M$ be a piecewise $C^{\infty}$ curve on $(M,F)$ with the velocity $\frac{d\gamma}{dt} = \frac{d{\gamma}^i}{dt} \frac{\partial}{\partial x^i} \in T_{\gamma (t)} M$. The arc length parameter of $\gamma$ is denoted by $s(t) = \int_{t_0}^{t} F(\gamma , \frac{d\gamma}{dr}) dr$, and the integral length  is given by $L(\gamma):= \int_a^b F(\gamma , \frac{d\gamma}{dt}) dt$.
 For every $x_0$ , $x_1$ $\in M$, denote by $\Gamma (x_0 , x_1)$ the collection of all piecewise $C^{\infty}$ curves  $\gamma : [a , b] \rightarrow M$ with $\gamma (a) = x_0$ and $\gamma(b) = x_1$, and define a map $d_F : M \times M \rightarrow [0 , \infty)$ by $d_F(x_0 , x_1) := inf L(\alpha)$, where $\alpha \in \Gamma (x_0 , x_1)$. It can be shown that $d_F$ satisfies the first two axioms of a metric space. Namely, (I) $d_F(x_0 , x_1) \geq 0$ , where equality holds if and only if $x_0 = x_1$; (II) $d_F(x_0 , x_1) \leq d_F(x_0 , x_1) + d_F(x_1 , x_2)$.\\
We should remark that the distance function $d_F$ on a Finsler space does not have  the symmetry property.
If the Finsler structure $F$ is reversible, we also have the third axiom of a metric space,
  (III) $d_F(x_0, x_1) = d_f(x_1,x_0)$, in this case $d_{F}$ is called \emph{symmetric}.
The \emph{Ricci scalar}  of $F$ is a positive scalar function given by  $Ric:=\frac{1}{F^2} R^i_i$. It is zero homogeneous  in $y$ which yields $Ric (x , y )$ depends on the direction of the flag pole $y$ but not to the length of $y$. The Ricci tensor of a Finsler
structure $F$ is defined by $Ric_{ij}:=\{\frac{1}{2}R^k_k\}_{y^iy^j}$, cf., \cite{Ak,BCS}.
A Finsler metric is called \emph{Einstein} if the Ricci scalar \emph{Ric}  is a function of $x$ alone, or equivalently if we have $Ric_{ij}=Ric(x)g_{ij}$.
\section{ A brief overview on projective changes of Finsler metrics}\label{sec;2+1}
In Finsler geometry, as discussed earlier the concept of a geodesic involves both the notion of a curve as a set of points and an orientation.
Two Finsler structures $F$ and $\bar{F}$ on a connected manifold $M$ are said to be  \emph{projectively related} if (considering the orientation) every geodesic of $F$ is a  geodesic of $\bar{F}$ and vice versa.\\
Let $\gamma: [a,b]\rightarrow M$ be an extremal of the length function (\ref{lenght}) and the length function
 $ L_{\bar F(c)}:= \int_a^b \bar F (c(t),\dot{c}(t)) dt$,  in the set of all smooth curves $c:=[a,b]\rightarrow M$ connecting $\gamma(a)$ and $\gamma(b)$, as well.
 Let $``s"$ be  the arc length parameter of $\gamma$ with respect to the $F$.  For $a\leq t \leq b$  we have $0\leq s(t) \leq A$ for some $A\in \mathbb{R}$. Now we consider  $\gamma: [0,A]\rightarrow M$ with parameter $``s"$. The arc length parameter, $``\bar {s}"$, of $\gamma$ with respect to the $\bar{F}$ is obtained by $\bar{s}=\int_0^s \bar{F}(x(r), \dot{x}(r)) dr.$
 Here, $\gamma$ is a  geodesic of $\bar{F}$, and we have
\begin{equation}\label{gpro1}
\frac{d^2{x}^i}{d\bar{s}^2}+2 \bar{G}^i(x(\bar{s}), \dot{x}(\bar{s}))=0.
\end{equation}
  On the other hand,  $$\frac{dx^i}{ds}=\frac{dx^i}{d\bar{s}} \frac{d\bar{s}}{ds},$$ $$\frac{d^2x^i}{ds^2}=\frac{d^2x^i}{d\bar{s}^2} (\frac{d\bar{s}}{ds})^2 + \frac{dx^i}{d\bar{s}} \frac{d^2 \bar{s}}{ds^2}.$$
 Applying the two latter equations to (\ref{geo}), we have
 $$\frac{d^2x^i}{d\bar{s}^2} (\frac{d\bar{s}}{ds})^2 + \frac{dx^i}{d\bar{s}} \frac{d^2 \bar{s}}{ds^2} + 2 G^i(x(s), \dot{x}(s))=0.$$
 Considering (\ref{gpro1}), we have
  \begin{equation}\label{gpro2}
2 G^i(x(s), \dot{x}(s))-2 \bar{G}^i(x(s), \dot{x}(s))=-\frac{\frac{d^2 \bar{s}}{ds^2}}{(\frac{d\bar{s}}{ds})^2} \frac{dx^i}{d\bar{s}}.
\end{equation}
It worths noting the fact that
\textit{``At each point of a connected manifold $M$, we have a precompact coordinate neighborhood $U$ such that given any $x \in U$ and $y \in T_x M$, there exists a unique geodesic $\sigma_{x,y} (t)$ which passes through $x$ at $t=0$ with the velocity $y$. Furthermore,  $\sigma_{x,y} (t)$ is $C^{\infty}$ in $``t"$ and in its initial data $x,y\neq 0$."}\cite{BCS}.\\
Considering this fact and the equations (\ref{gpro2}), it is obvious that if two Finsler metrics are projectively invariant, there exists  a  1-homogeneous scalar field $P(x,y)$ satisfying
\begin{equation}\label{pro*}
\bar{G}^i(x,y)=G^i(x,y)+ P(x,y)y^i.
\end{equation}
Now, let consider that there exists  a  1-homogeneous scalar field $P(x,y)$ satisfies (\ref{pro*}).
 Let $\gamma:= (x^i(t))$ be a geodesic of the Finsler metric $F$ on $M$. We take parameter ``$\bar {t}$" determined by
 $$\frac{d^2\bar{t}}{dt^2}=2p(x^i(t),\dot{x}(t))\frac{d \bar{t}}{dt}.$$By direct calculation, we have $$\frac{d^2x^i}{d{\bar{t}}^2}+2{\bar{ G}}^i(x(\bar{t}), \dot{x}(\bar{t}))=0.$$
Therefore, $\gamma$ is a geodesic of $\bar{F}$ as well.
Considering the above arguments,\\ \textit{Two Finsler metrics $F$ and $\bar{F}$ are projectively equivalent if and only if there exists  a  1-homogeneous scalar field $P(x,y)$ satisfies
 (\ref{pro*}).}\\
\\
\section{Schwarzian derivative and a pseudo-distance}\label{sec;3}
    Let $\gamma=x^i(t)$ be a geodesic on $(M,F)$. In general, the parameter ``$t$" of $\gamma$, does not remain invariant under projective changes of $F$.
     A  parameter  is called  \emph{projective}, if it remains invariant under a projective change of $F$.
In \cite{B,E,T} the projective parameters are defined  for geodesics of the general affine connections. In \cite{SB2}
 the projective parameters are carefully spelled out for geodesics of a Finsler metric  as a solution to the following ODE
$$\{p,s\}:= \frac{\frac{d^3p}{ds^3}}{\frac{dp}{ds}}-\frac{3}{2}\big[\frac{\frac{d^2p}{ds^2}}{\frac{dp}{ds}}\big]^2 = \frac{2}{n-1} Ric_{jk}\frac{d{x}^j}{ds}\frac{d{x}^k}{ds},$$
where $\{p,s\}$ is  the well known \emph{Schwarzian derivative}
 and  $s$ is the arc length parameter of $\gamma$.
 The projective parameters are unique up to a  linear fractional transformation
$
 \{ \frac{a p + b}{c p +d} , s \} = \{p , s\},$
where, $ad-bc\neq 0$.\\Here, we briefly introduce the pseudo-distance $d_M$ and point out some of its properties. See  \cite{SB2} for all details.
Consider the open interval  $I=(-1,1)$ with Poincar\'e metric $ds^2_I=\frac{4du^2}{(1-u^2)^2}$. The Poincar\'e distance between two points $a$ and $b$ in $I$ is given by
$\rho(a , b) = \mid \ln \frac{(1-a)(1+b)}{(1-b)(1+a)} \mid,$
cf.,  \cite{O}.
 A geodesic $f :I \rightarrow M$ on the Finsler space $(M,F)$ is said to be \emph{projective}, if the natural parameter $u$ on $I$ is a projective parameter.

Given any two points $x$ and $y$ in $(M,F)$, we consider a chain $\alpha$ of geodesic segments joining these points. That is; i) a chain of points $x = x_0 , x_1 , ... ,x_k = y$ on $M$; ii)
 pairs of points $a_1,b_1 ,..., a_k,b_k$ in $I$; iii)
 projective maps $f_1,...,f_k$, $f_i: I \rightarrow M $, such that
$f_i(a_i) = x_{i-1}, \quad f_i(b_i) = x_i, \quad i = 1,...,k.$
By virtue of the Poincar\'{e} distance $\rho(.,.)$ on $I$ we define the length $L(\alpha)$ of the chain $\alpha$ by
$L(\alpha):= \Sigma_i \rho(a_i , b_i)$, and we put $d_M(x , y):= inf L(\alpha)$, where the infimum is taken over all chains $\alpha$ of geodesic segments from $x$ to $y$.

Based on a projectively invariant pseudo-distance defined by the Schwarzian derivative  on a Finsler manifold, the   Ricci scalar   is studied under a projective change of $F$.
 Next, the Schwarzian derivative  is applied to introduce  a projectively invariant parameter, called \textit{projective parameter}, for a given geodesic on $(M,F)$.
 Then a certain chain of piecewise smooth geodesics  $\alpha$, which joins a point $x$ in $M$ to a point $y$ in $M$ is considered.
 The \textit{Poincar\'e} metric on the open interval $I$ determines $L(\alpha)$  the length of $\alpha$. Finally, the projectively invariant pseudo-distance $d_M$ is defined to be the infimum of $L(\alpha)$ for different chains $\alpha$, see \cite{SB2,SB1} for more details.
 Let us review the following properties
\setcounter{thm}{0}
\begin{Aproposition}\label{Pro;1}\cite{SB2}
Let $(M, F)$ be a Finsler space. Then for any points $x$, $y$, and $z$ in $M$,  $d_M$ satisfies
 \begin{itemize}
 \item[(i)]$d_M(x,y)=d_M(y,x)$.
 \item[(ii)]$d_M(x , z) \leq d_M(x , y) + d_M(y , z)$.
 \item[(iii)]If $x = y$ then $d_M(x , y) = 0$ but the inverse is not  always true.
 \end{itemize}
\end{Aproposition}

\begin{Aproposition}\label{Prop;SHI}\cite{Shi}
Let $(M, F)$ be  a connected, complete Finsler manifold with associated distance function $d_F:M \times M \to [0,1)$.
Then $d_F$ is symmetric distance function on $M \times M$ if and only if $F$ is
absolute homogeneous, that is, $F(x, y)= F(x,-y)$.
\end{Aproposition}
\begin{Atheorem}
\cite{SB2} Let $(M,F)$ be a connected Finsler manifold where the Ricci tensor satisfies
\begin{equation*}
(Ric)_{ij} \leq -c^2g_{ij},
\end{equation*}
as matrices, for a positive constant $c$. Then the pseudo-distance $d_M$, is a distance.
\end{Atheorem}
\begin{Atheorem}\cite{SB3} Let $(M,F)$ be a connected complete Finsler space of positive semi-definite  Ricci tensor. Then the intrinsic projectively invariant pseudo-distance $d_M$ is trivial, that is $d_M=0$.
\end{Atheorem}
\begin{Atheorem}\label{TH;C} \cite{SB3}
 Let $(M,F)$ be a connected (complete) Finsler space. If the  Ricci tensor  is negative-definite and parallel with respect to the Berwald or Chern connection, then  the intrinsic projectively invariant pseudo-distance $d_M$, is a (complete) distance.
\end{Atheorem}
\begin{Atheorem}\cite{SB1}\label{THD}
 Let $(M , F)$ be a complete Einstein Finsler space with
\begin{equation*}
(Ric)_{ij} = -c^2 g_{ij},
\end{equation*}
 where $c$ is a positive constant. Then  the projectively invariant distance $d_M$, is proportional to the Finslerian distance $d_F$, that is
 \begin{equation*}
 d_M(x , y) = \frac{2c}{\sqrt{n - 1}} d_F(x , y).
 \end{equation*}
 \end{Atheorem}
\section{Reversibility and Ricci curvature}
In this section we study some necessary conditions for reversibility of a Finsler structure. We give a proof for Theorem \ref{Th;main} and Theorem \ref{Th;2} and investigate global conclusions.
\\

\textbf{Proof of Theorem\ref{Th;main}.}
Let $(M,F)$ be a connected  complete Finsler space of constant negative Ricci scalar.
According to part $(i)$ of Proposition\ref{Pro;1}, $d_M$ is symmetric.  Apply this property to Theorem \ref{THD} reads that $d_F$ is symmetric  too.  After Proposition \ref{Prop;SHI} on a  complete connected Finsler manifold the associated distance function $d_F$ is symmetric if and only if $F$ is reversible.
The latter fact completes the proof.
\hspace{\stretch{1}}$\Box$\\
\bc\label{Th;1}
 Let $(M,F)$ be a complete Finsler manifold of constant negative flag curvature. Then the Finsler structure $F$ is reversible.
\ec

The assumption of reversibility of a Finsler structure excludes many interesting examples such as Randers metrics. A Randers metric on a manifold $M$ is a Finsler structure in the form $F=\alpha + \beta$,
where $\alpha = \sqrt{a_{ij}(x)y^iy^j}$ is a Riemannian metric and $\beta = b_i(x)y^i$ is a one form on $M$, where $b:= || \beta_{x} ||-{\alpha} <1$. It is well-known that a Randers metric is reversible if and only if it is Riemannian. Given the above statement, we point out that a non-Riemannian Randers metric can only adopt one of the properties: completeness and having a constant negative Ricci scalar.
Therefore we have Corollary\ref{Cor;1}.

%
Let us $\bigtriangledown^b$ denotes the covariant derivative of Berwald connection and the Ricci tensor be parallel with respect to the Berwald connection. The following  similar arguments hold as well for Chern connection. We have
 \begin{equation}\label{Ber1}
\bigtriangledown^b_{\frac{\delta}{\delta x^j}}Ric_{hl}=\frac{\delta Ric_{hl}}{\delta x^j}- Ric_{hr}{G^r}_{lj}-Ric_{lr}{G^r}_{hj}=0,\quad {G^r}_{lj}=1/2\frac{\partial^2 G^r}{\partial y^j y^l}.
\ee
\begin{equation}\label{Ber2}
\bigtriangledown^b_{\frac{\partial}{\partial y^k}}Ric_{ij}=\frac{\partial Ric_{ij}}{\partial y^k}=0.
\ee
Contracting (\ref{Ber1}) in $Ric^{ih}y^jy^l$ we have
\begin{equation}\label{Ber3}
Ric^{ih}y^jy^l\frac{\partial Ric_{hl}}{\partial x^j}-G^i-1/2Ric^{ih} Ric_{la}\frac{\partial G^a}{\partial y^h}y^l=0.
\ee
On the other hand
$$-1/2 Ric^{ih}y^ry^l\frac{\partial Ric_{lr}}{\partial x^h}+ 1/2Ric^{ih}y^ry^lRic_{la}{G^a}_{rh}+1/2 Ric^{ih}y^ry^lRic_{ra}{G^a}_{lh}=0.$$
\begin{equation}\label{Ber4}
-1/2 Ric^{ih}y^ry^l\frac{\partial Ric_{lr}}{\partial x^h}+ 1/2Ric^{ih}y^r Ric_{ra}\frac{\partial G^a}{\partial y^h}=0.
\ee
Considering (\ref{Ber3}) and (\ref{Ber4}) we have ${\hat{G}}^i=G^i$. Therefor  $Ric_{ij}$ is the Ricci tensor of $\hat{F}$ too.

The above discussions are spelled out in the following result.
\setcounter{thm}{1}
\bt\label{Th;2}
 Let $(M,F)$ be a complete connected  Finsler space. If the Ricci tensor is negative-definite
and parallel with respect to Berwald or Chern connection then the geodesics are reversible.
\et
 \bpf
 Adding completeness to the  hypothesis of Theorem \ref{TH;C}  and applying Theorem \ref{Th;main}, we can prove Theorem \ref{Th;2}.  Put $\hat{F}(x,y)=\sqrt{-Ric_{ij}(x,y)y^iy^j}$.  We have
  $\hat{g}_{ij} = {[1/2\hat{F}^2]}_{y^iy^j}= -Ric_{ij}(x,y)$.


 Moreover,
\be\label{affin}
\hat{g}_{ij} = - Ric_{ij}=- \hat{Ric}_{ij}.
\ee
Applying (\ref{affin})  to Theorem \ref{Th;main}, reads $\hat{F}$ is reversible and hence it is geodesically reversible. The Finsler manifolds $(M,F)$ and  $(M,\hat{F})$ are affine, therefore the Finsler structure $F$ is geodesically reversible.
\hspace{\stretch{1}}
\epf
\br
The above discussion also provides the opportunity to characterize the  completeness of a large class of Finsler metrics with a constant negative Ricci scalar.
For instance we know that the Finslerian Pioncar\'e metric is of constant flag curvature  $-1/4$.
 Therefore, by the above discussion, it could not be complete.
 Forward and backward completeness of this metric is discussed in details in \cite{BCS}.
 \er
 Meanwhile, the above discussion gives a characterization of the Einstein Finsler spaces as follows
\br\label{Rem;2}
It  is noteworthy to mention that the obtained result on Corrolary\ref{Cor;1} does not contradict the examples given by X. Cheng and Z. Shen in  \cite{chsh}, where they have shown; there are singular Einstein $(\alpha, \beta)$-metrics with $Ric=0, Ric=-1$ and $Ric=+1$ on the sphere $S^3$.

 In fact, as the authors have explicitly explained at the end of the \cite{chsh},  the results are held on Finsler spaces with \emph{almost regular}  metrics and fail to be true in the case of \emph{regular} metrics.
\er
\textbf{ Acknowledgement.}
The first author would like to thank the Institute of Mathematics Toulouse France (ITM) where this article is partially written.

  Behroz Bidabad\\
Professor of Department of Mathematics and Computer Sciences\\
Amirkabir University of Technology (Tehran Polytechnic),
424 Hafez Ave. 15914 Tehran, Iran.
E-mail: bidabad@aut.ac.ir\\
\vspace{0.5cm}

 Maryam Sepasi\\
Faculty of Science, Islamic Azad University of Shiraz
Sadra, Shiraz 71993, Iran.
E-mail: sepasi@iaushiraz.ac.ir


\begin{thebibliography}{}
 \bibitem{Ak}   Akbar-Zadeh,  H.; Sur les espaces de Finsler \`{a} courbures sectionnelles
 constantes ,  Acad .  Roy .  Belg .  Bull .  CI .  Sci .  (5) 14 ,  (1988) ,  281-322 .
  \bibitem{B}  Berwald, L.; On the projective geometry of paths ,  Proc .  Edinburgh Math .  Soc .  (5) ,  (1937) ,   103-115.
   \bibitem{BCS}  Bao, D.,  Chern, S.S.,    Shen, Z.;   Riemann-Finsler geometry ,
 Springer-Verlag ,  2000.
  \bibitem{chsh} Cheng, X. ,  Shen, Z.; Einstein Finsler metrics and Killing vector fields on Riemannian manifolds, Sci. China Math. 60 (2017), no. 1, 83–-98.
 \bibitem{SB2} Bidabad, B. ,    Sepasi, M.; On a projectively invariant pseudo-distance in Finsler geometry, Int. J. Geom. Methods in Modern Phys. Vol. 12, No. 04, (2015), 1550043--1550055.
\bibitem{SB3}  Bidabad, B. ,    Sepasi, M.; On the characteristic of projectively invariant pseudo-distance on Finsler spaces, IJST, 39, A2 (2015) , 233-238.
 \bibitem{E}  Eisenhart, L.P.; Non-Riemannian geometry, AMS Colloquium Publications, Volume 8, 1927.
  \bibitem{FF}Fang,Y. Foulon, P.; On Finsler manifolds of negative flag curvature. J. Topol. Anal. 7 (2015), no. 3, 483–-504.
 \bibitem{F1} Foulon, P.; Locally symmetric Finsler spaces in negative curvature, C. R. Acad.
Sci. Paris, ser. I Math., no.10, 324 (1997), 1127-1132.
 \bibitem{F2} Foulon, P.; Curvature and global rigidity in Finsler manifolds, Houston J. Math.,
28, (2002), 263-292.
\bibitem{Ma} Matveev, V.; On projective equivalence and pointwise projective relation of Randers metrics, Internat. J. Math. 23 (2012), no. 9, 1250093, 14 pp.
\bibitem{MSh} Mo, X.,  Shen, Z.; On negatively curved Finsler manifolds of scalar curvature, Canadian Mathematical Bulletin, 48, (2005), 112-120.
 \bibitem{O} Okada, T.; On models of projectively flat Finsler spaces of constant negative curvature," Tensor, N. S. Vol 40 (1983) 117-124.
\bibitem{SB1} Sepasi, M., Bidabad, B.; On a projectively invariant distance on Finsler spaces of
 constant negative Ricci scalar, C. R. Acad. Sci. Paris, Ser. I. Vol. 352, Issue 12, ( 2014), 999-1003.
\bibitem{Shi} Sabau, Sorin V., and Shimada, H.; Finsler manifolds with reversible geodesics, Rev. Roumaine Math. Pures Appl. 57 (2012), no. 1, 91-–103.
  \bibitem{Sh} Shen, Z.; Finsler manifolds with non-positive flag curvature and constant S-curvature, Math. Z. Volume no. 3, 249 (2005), 625-639.
\bibitem{Sh2} Shen, Z.; On projectively related Einstein metrics in Riemann-Finsler geometry,  Math. Ann. 320,625--647 (2001).

          \bibitem{T} Thomas, T.Y.; On the projective and equi-projective geometries of paths, Proc. N. A . S. 11 (1925).

 \end{thebibliography}
\end{document}